\documentclass[11pt]{amsart}
\usepackage{amsmath,amsfonts,amsthm,amssymb,eufrak,comment}
\usepackage{amscd}
\usepackage{times,fancyhdr}
\usepackage{array}
\usepackage{epsfig}
\usepackage{graphicx}
\usepackage{color}
\usepackage{mathtools}
\usepackage{mathabx}

\usepackage{tikz-cd}

\usepackage{graphicx}
\usepackage{tikz}
\usetikzlibrary{decorations.fractals}
\usetikzlibrary{lindenmayersystems}
\usepackage{pgfplots}

\usepackage[matrix,arrow,ps]{xy}
 \UsePSspecials{dvips}
\usetikzlibrary{positioning}

\input{xy}

\usepackage{euscript}

\usepackage{footnote}
\makesavenoteenv{tabular}
\makesavenoteenv{table}

\newcounter{thecounter}
\numberwithin{thecounter}{section}

\newtheorem{proposition}[thecounter]{Proposition}
\newtheorem{theorem}[thecounter]{Theorem}

\DeclareMathOperator{\GL}{GL}

\DeclareMathOperator{\Aut}{Aut}

\DeclareMathOperator{\End}{End}
\DeclareMathOperator{\Exp}{Exp}

\DeclareMathOperator{\Ad}{Ad}
\DeclareMathOperator{\ad}{ad}

\DeclareMathOperator{\re}{{re}}
\DeclareMathOperator{\im}{{im}}

\renewcommand{\a}{\alpha}

\newcommand{\Z}{{\mathbb Z}}

\newcommand{\C}{{\mathbb C}}
\newcommand{\N}{{\mathbb N}}

\newcommand{\RomanNumeralCaps}[1]
    {\MakeUppercase{\romannumeral #1}}

\begin{document}

\title{\bf{Constructing a Lie group analog for \\the Monster Lie algebra}}

\date{\today}

\begin{abstract} Let $\frak m$ be the  Monster Lie algebra. We summarize several interrelated constructions of Lie group analogs for  $\frak m$. Our constructions are analogs for $\frak m$ of Chevalley and Kac--Moody groups and their generators and relations.
\end{abstract}

\author[Lisa Carbone]{Lisa Carbone}
\address{Department of Mathematics, Rutgers University, Piscataway, NJ 08854-8019, USA}
\email{lisa.carbone@rutgers.edu\footnote{Corresponding author}}

\author[Liz Jurisich]{Elizabeth Jurisich}
\address{Department of Mathematics, College of Charleston,
66 George Street,
Charleston, SC 29424, USA}
\email{jurisiche@cofc.edu}

\author[Scott H. Murray]{Scott H. Murray}
\address{University of Toronto Mississauga,
3359 Mississauga Road, Mississauga, ON L5L 1C6, Canada}
\email{scotthmurray@gmail.com}

\thanks{The first author's research is partially supported by the Simons Foundation, Mathematics and Physical Sciences-Collaboration Grants for Mathematicians, Award Number 422182}
\maketitle


\section{Introduction}
Let $\frak m$ be the Monster Lie algebra,   constructed by Borcherds (\cite{B1, B2, B3}). Borcherds discovered $\frak m$ as an example of a  new class of Lie algebras which became known as {\it generalized Kac--Moody algebras}, or {\it Borcherds algebras}. The  invariant bilinear form on $\frak m$ gives rise to the following generalized Cartan matrix  (\cite{Jur1}):
\begin{center}
  \includegraphics[scale=0.42]{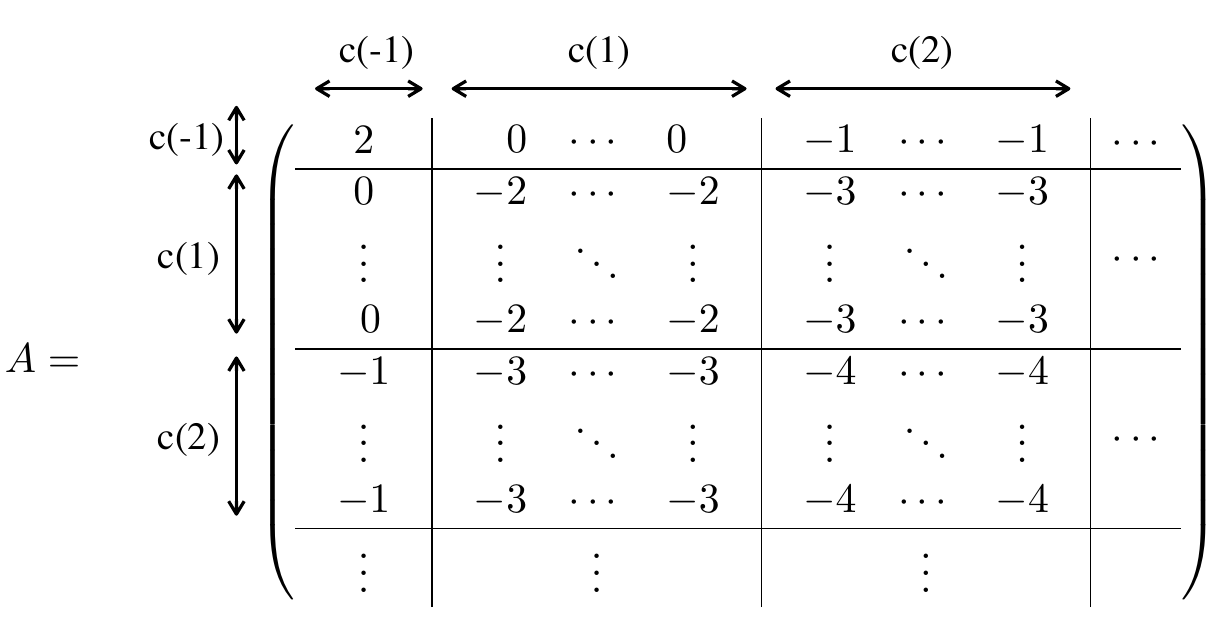}
\end{center}
where  $c(i)$ is the coefficient of $q^i$ 
in the  modular function
$$J(q)=j(q) -744=\sum_{i\geq -1}c(i)q^i= {1 \over q}  + 196884 q + 21493760 q^2 + 864299970 q^3  + \cdots,$$
and $j(q)$ is the well-known $j$-function.

 Borcherds constructed $\frak m$  (\cite{B2})  to prove part of the Conway--Norton conjectures (\cite{CN}). A fundamental component of Borcherds' construction was Frenkel, Lepowsky and Meurman's representation (\cite{FLM1}, \cite{FLM2}) of the Fischer--Griess Monster  finite simple group {\bf M}.
Borcherds constructed $\frak m$ as a certain  quotient of the `physical space'  of the vertex operator algebra $V=V^\natural\otimes V_{{1,1}}$, where $V^\natural$ the moonshine module of \cite{FLM1}, \cite{FLM2},   a graded $\bold M$-module with $\Aut(V^\natural)=\bold M$, and 
 $V_{{1,1}}$ is  a vertex operator algebra  for the even unimodular  2-dimensional Lorentzian lattice  $\RomanNumeralCaps{2}_{1,1}$.

The Monster Lie algebra has an equivalent characterization as the Lie algebra $\mathfrak m =\mathfrak g(A) / \mathfrak z$ given by generators and relations associated to the matrix $A$, where $\frak z$ is the center of $\frak g(A)$  (\cite{B3}, \cite{Jur1} and Section~\ref{gensrelns}).  We use this latter characterization of $\frak m$ for our purposes.

The Monster Lie algebra  has the  triangular decomposition $\frak m= \mathfrak{n}^-\oplus \mathfrak{h}\oplus \mathfrak{n}^+$. 
We also use the following  decomposition 
$$\frak m\cong \mathfrak{u}^-\oplus \mathfrak{gl}_2(-1)\oplus \mathfrak{u}^+$$
where $\mathfrak{gl}_2(-1)\cong \mathfrak{gl}_2$, and $\mathfrak{u}^+$ (resp.\ $\mathfrak{u}^-$) are subalgebras freely generated by certain positive (respectively negative) imaginary root vectors (\cite{Jur1}, \cite{Jur2} and Theorem~\ref{JurDecomp}).

 The  class of Borcherds algebras has been widely studied. For example, the appearance of $\mathfrak{m}$ and other Borcherds algebras as symmetries in heterotic string theory have been noted (\cite{Ca}, \cite{HM}, \cite{PPV}). However, there have been no constructions of Lie group analogs for Borcherds algebras.

There are many methods for constructing Kac--Moody groups, the Lie group analogs for infinite dimensional Kac--Moody algebras (see, for example, \cite{CG, Ga1, GW, KP, Ku, Ma, Mar, Ma84, MT, Rou, Sl, Ti}).  For Borcherds algebras, the situation is quite different, due in part to the absence of suitable representation theoretic, geometric, algebraic geometric and analytic methods in this setting.

In particular, the  root vectors associated with real roots of a Kac--Moody or Borcherds algebra are locally $\ad$-nilpotent. In the Kac--Moody case,  the  simple roots are all real and associated choices of root vectors generate the whole Kac--Moody algebra. For Borcherds algebras, this is no longer the case. There are standard Chevalley generators of $\mathfrak{n}^\pm$ associated to simple imaginary  roots, but these  do not act locally nilpotently on the adjoint representation, or on the highest weight modules of interest.

Here we report on several works in progress where we  construct groups associated to $\frak m$. We give several interrelated constructions which serve different purposes. 

We  construct a complete pro-unipotent group of automorphisms of a completion $\widehat{\frak m}$ of $\frak m$ (Section~\ref{prounipgp} and \cite{CJM}) by proving convergence of these maps on finite dimensional subspaces. For this purpose, we introduce the notion of \emph{pro-summability} (Section~\ref{prosum} and \cite{CJM}). 

We also construct a group $G(\frak m)$ given by generators and relations, and show that 
it has an analog of a unipotent subgroup which acts as automorphisms of  $\widehat{\frak m}$ (Section~\ref{gensrelnsgp} and \cite{CJM}, \cite{ACJM}).

To construct the analog of a simply connected  Kac--Moody Chevalley group for $\frak m$  would require the use of an integrable highest weight representation of $\mathfrak m$. Since the simple imaginary root vectors which generate  $\frak m$ are not locally nilpotent on any of the highest weight modules that we encounter (including generalized Verma modules), we need a different approach. It was shown in \cite{JLW} that  $\frak m$  has a representation on a certain  tensor algebra $T(\mathcal V)$ analogous to  standard irreducible modules  for semisimple and Kac--Moody algebras. The  parabolic  subalgebra $ \frak{gl}_2(-1)\oplus \frak u^+$ of $\frak m$ is locally nilpotent on the $\frak m$-module $T(\mathcal V)$ (\cite{JLW}). Using this representation, we construct the analog of a simply connected  Kac--Moody Chevalley group associated to the parabolic subalgebra $ \frak{gl}_2(-1)\oplus\frak u^+$ of $\frak m$ (Section~\ref{simplyconn} and \cite{CGJM}). 


The authors would like to thank Abid Ali, Alejandro Ginory, Yi-Zhi Huang, James Lepowsky, Ugo Moschella and Siddhartha Sahi for their interest in this project and for helpful discussions. Much of the early part of this work was undertaken at Universit\'a degli studi dell'Insubria. The first two authors gratefully acknowledge the hospitality and discussions with Ugo Moschella that made this work possible. We would also like to express our gratitude to the referee for careful reading of the paper.

\section{The Monster Lie algebra}\label{gensrelns}


For $j\in\Z$, we recall the definition of  $c(j)$ above. Define index sets 
$$I^{\re}=\{(-1,1)\}, \quad I^{\im} = \{(j,k)\mid j\in\N,\; 1 \leq k\leq c(j)\},$$
$$I = \{(j,k)\mid j,k\in\Z,\; 1 \leq k\leq c(j)\} =I^{\re} \sqcup I^{\im}$$
and generalized Cartan matrix
$$ A = ( a_{jk,pq} )_{(j,k),(p,q)\in {I}}, \quad \text{where $a_{jk,pq} = -(j+p)$}.$$  
The Lie algebra $\mathfrak{g}(A)$  has generating set
 $\{{e}_{jk},\ {f}_{jk},\ h_{jk}\mid (j,k)\in {I}\} $
and defining relations 
\begin{align*}
 [h_{jk},h_{pq}]&=0,\\
 [h_{jk},{e}_{pq}]&=-(j+p) {e}_{pq},\\
 [h_{jk},{f}_{pq}]&=(j+p) {f}_{pq},\\
 [{e}_{jk},{f}_{pq}]&=\delta_{jp}\delta_{kq}h_{jk},\\ 
 (\ad {e}_{-1\,1})^j \,{e}_{jk}&= (\ad {f}_{-1\,1})^j \,{f}_{jk}=0,
\end{align*}
for all $(j,k),\,(p,q) \in {I}$. 
The generators  are  the usual Chevalley generators of the Lie algebra associated with the matrix $A$, where
the double index is chosen to reflect the block structure of  $A$ (see also \cite{JLW}, page 10). The Cartan subalgebra is  ${\mathfrak h}_A := \sum_{(j,k)\in{I}} \C h_{jk}\subseteq \frak g(A)$.

As in \cite{Jur1}, \cite{Jur2}, we define the Monster Lie algebra $\mathfrak{m}$  to be  $\frak m := \frak g(A)/\frak z$ where $\frak z$ is the center of the Lie algebra.
Define the following elements of $\frak m$:
\begin{align*}
h_1 &:= \frac12(h_{-1\,1}-h_{1\,1})+\mathfrak{z}, & h_2 &:= -\frac12(h_{-1\,1}+h_{1\,1})+\mathfrak{z},\\
e_{-1} &:= {e}_{-1\,1}+\mathfrak{z}, & f_{-1} &:= {f}_{-1\,1}+ \mathfrak{z}.
\end{align*}
We will write $e_{jk}$ for $e_{jk}+\mathfrak{z}$ and $f_{jk}$ for ${f}_{jk}+ \mathfrak{z}$, for  all $(j,k)\in {I}^{\im}$.

Then $\frak m$ is the 
Borcherds (generalized Kac--Moody) algebra with generating set 
$$\{e_{-1},f_{-1},h_{1},h_2\}\cup \{e_{jk},\ f_{jk}\mid (j,k)\in I^{\im}\}.$$
These generators are subject to the following defining relations (\cite{Jur1}, \cite{Jur2}):
\begin{align*}
[h_{1},h_2]&=0,\\
[h_1,e_{-1}] &= e_{-1},&              
              [h_2, e_{-1}] &= -e_{-1},\\
 [h_1,e_{jk}]&= e_{jk},&
              [h_2, e_{jk}] &= j e_{jk},\\
[h_1,f_{-1}] &= -f_{-1},&
             [h_2,f_{-1}] &= f_{-1},\\
 [h_1,f_{jk}]&= - f_{jk},&
              [h_2,f_{jk}] &= -j f_{jk},\\
 [e_{-1},f_{-1}]&=h_1-h_2, \\
 [e_{-1},f_{jk}]&=0,  &[e_{jk},f_{-1}]&=0,\\
 [e_{jk},f_{pq}]&=-\delta_{jp}\delta_{kq} \left(jh_1 + h_2\right),\\ 
(\ad e_{-1})^j e_{jk}&=0,&\qquad (\ad f_{-1})^j f_{jk}&=0,
\end{align*}
for all $(j,k),\,(p,q) \in {I}^{\im}$. The Cartan subalgebra of $\frak m$ is ${\frak h} := {\frak h_A}/{\frak z}=\C h_1\oplus \C h_2\subseteq \frak m$. Note that $A$ has rank $2$. 

Now $e_{-1}$ (resp.\ $f_{-1}$) is the positive (respectively negative) simple real root vector,
and the generators $e_{jk}$ (resp.\ $f_{jk}$) are the positive (respectively negative) simple imaginary root vectors.

We define the extended index set
$$ E = \{(\ell,j,k)\mid( j,k)\in I^{\im},\, 0\leq \ell<j \} = \{(\ell,j,k)\mid j\in\N,\, 1 \leq k\leq c(j),\, 0\leq \ell<j \}.$$
and set
$$
   e_{\ell,jk}:= \frac{(\ad e_{-1})^\ell e_{jk}}{\ell!}\quad\text{and}\quad f_{\ell,jk}:= \frac{(\ad f_{-1})^\ell f_{jk}}{\ell!},
$$
for $ (\ell,j,k) \in E$.

The following non-trivial result gives an additional non-standard generating set which we will use for our group constructions.
\begin{theorem}\label{JurDecomp}(\cite{Jur1}, \cite{Jur2}) Let $\mathfrak{gl}_{2}(-1)$ be the subalgebra of $\frak m$ with basis $\{e_{-1}, f_{-1}, h_{1}, h_2\}$. Then
$$\mathfrak{m} = \mathfrak{u}^- \oplus \mathfrak{gl}_{2}(-1)  \oplus \mathfrak{u}^+$$
where
$\mathfrak{gl}_{2}(-1):= \langle e_{-1},f_{-1},h_1,h_2\rangle\cong \mathfrak{gl}_2$,  $\;\mathfrak{u}^+$ is a subalgebra freely generated by  $\{e_{\ell,jk}\mid (\ell,j,k)\in E\}$  and
$\mathfrak{u}^-$ is a subalgebra freely generated by  $\{f_{\ell,jk}\mid (\ell,j,k)\in E\}$.
\end{theorem}




\section{$\exp$ and $\ad$ for infinite dimensional Lie algebras}\label{expad}

For a finite dimensional semisimple Lie algebra $\frak g$, the simple root vectors $e_i$ and $f_i$ act nilpotently on the adjoint 
representation, that is,  $e_i$ and $f_i$ are $\ad$-nilpotent. Thus there are automorphisms of the form $$\exp(u\,\ad(x)) (y)=y +u[x,y]+\dfrac{u^2}{2!}[x,[x,y]]+\cdots$$
For $x=e_i$ or $x=f_i$, these are finite sums, thus  are well defined automorphisms of $\frak g$.

For a Kac--Moody algebra $\mathfrak{g}$, the simple root vectors $e_i$ and $f_i$ are real and hence they act {\it locally nilpotently} on
$\frak{g}$. That is,  for all $y\in  \mathfrak{g}$, we have $(\ad(e_i))^{n}(y)=0$ and  $(\ad(f_i))^{m}(y)=0$ for some $n,m\gg 0$. 
This means that $\exp(u\,\ad(e_i))$ and $\exp(v\,\ad(f_i))$, $u,v\in\C$, are summable. Recall that an infinite sum $\sum_n E_n$ of operators is called \emph{summable} if $\sum_n E_n(y)$ reduces to a finite sum for all $y$ (\cite{LL}).  Again this means that $\exp(u\,\ad(e_i))$ and $\exp(v\,\ad(f_i))$
can be viewed as  elements of $\Aut(\mathfrak{g})$.
Replacing $\ad$ with a representation on an integrable highest weight module $V$, a similar method gives an analog of a simply connected Chevalley group in the Kac--Moody case (\cite{CG}, \cite{CLM}).

For the Monster Lie algebra $\mathfrak{m}$, the real simple root vectors $e_{-1}$ and $f_{-1}$ act locally nilpotently on the full adjoint 
representation. However,  the imaginary simple root vectors $e_{jk}$ and $f_{jk}$ do not necessarily act locally $\ad$-nilpotently, and so
this approach no longer works for the Monster Lie algebra (or other Borcherds algebras with imaginary simple roots). 
The same is true when we replace the adjoint 
representation by  faithful highest weight modules, 
including generalized Verma modules.


\section{Completion $\widehat{\frak m}$ of $\frak m$ and pro-summability}\label{prosum}
The algebra $\frak m$  has the usual triangular decomposition
$$\frak m = \frak n^-\oplus \frak h \oplus\frak n^+$$
where  
$\frak n^{\pm} =\bigoplus_{\alpha\in \Delta^{\pm}}\frak m_{\alpha}$, $\frak h$ is the Cartan subalgebra, $\Delta^\pm$ are the sets of positive (respectively negative) roots and  $\frak m_\a$ are the root spaces. Let $\Delta=\Delta^+\sqcup\Delta^-$ and let $Q$ denote the root lattice, that is, the $\Z$-span of $\Delta$.

The algebra $\frak m$ is $Q$-graded.  We can define a $\Z$-grading on $\frak m$ by 
defining a map 
\begin{align*}
\lambda:Q &\to\Z,\\
(m, n) &\mapsto m + n.
\end{align*}
This gives compatible $\Z$-grading whose homogeneous components are of the form
$$\mathfrak{m}_k = \bigoplus_{\a\in\lambda^{-1}(k)\cap\Delta} {\mathfrak m}_{\a},$$
for $k\ne0$, and $\mathfrak{m}_0=\mathfrak{gl}_{2}(-1)$.
The set $\Delta_k:= \lambda^{-1}(k)\cap\Delta$ is finite, which ensures that $\mathfrak{m}_k$ is finite dimensional for each $k\in\Z$.
For $i>0$, let $$\mathfrak{n}_i=  \bigoplus_{i'\ge i} \mathfrak{m}_{i'}.$$
 We have a descending chain of ideals (\cite{CJM})
$$\mathfrak{n}^+ = \mathfrak{n}_0 \ge \mathfrak{n}_1 \ge \cdots\ge \mathfrak{n}_i\ge \cdots $$
with the property that  $\mathfrak{n}^+/\mathfrak{n}_i$ is a finite dimensional nilpotent Lie algebra for each $i\geq 0$. Thus $\frak n^+$ is  {\it pro-nilpotent}.

As in \cite{Ku}, Chapter IV, Section 4, 
we define the \emph{(positive) formal completion} of $\frak m$ to be
$$\widehat{\frak m} = \frak n^- \oplus \frak h \oplus \widehat{\frak n}^+,$$
where
$$\widehat{\frak n}^+ := \prod_{\a\in\Delta^+} \frak m_\a = \prod_{i\in\mathbb N} \frak m_i $$
is the pro-nilpotent completion of $\mathfrak{n}^+$. 
Let $$\widehat{\mathfrak{n}}_i=\prod_{i'\ge i} \mathfrak{m}_{i'}.$$ 
Then $$\widehat{\frak n}^+\cong \varprojlim_{i\geq 1}\widehat{\mathfrak{n}}^+/\widehat{\mathfrak{n}}_i.$$

Let  $e_{\ell,jk}$ be the imaginary  root vector as above. Then $\exp(u\,\ad(e_{\ell,jk}))(y)$ is generally an infinite sum for $u\in\C$, $y\in\frak m$. 
However, we call  $\exp(u\,\ad(e_{\ell,jk}))$  \emph{pro-summable}
since $\exp(u\,\ad(e_{\ell,jk})) (y)$ reduces to a finite sum for all $y$ in
$$\frak n^- \oplus \frak h \oplus (\widehat{\mathfrak n}^+/\widehat{\mathfrak{n}}_i)$$ 
and for all $i>0$
(\cite{CJM}).
By taking the inverse limit, $\exp(u\,\ad(e_{\ell,jk}))$ is a well-defined automorphism of $\widehat{\frak m}$.

Note, however, that $\exp(u\,\ad(e_{\ell,jk}))$ is \emph{not} a well-defined automorphism of ${\frak m}$. Nor can $\exp(u\,\ad(f_{\ell,jk}))$
be defined as an automorphism of $\frak m$ or $\widehat{\frak m}$, although it can be defined as an automorphism of the corresponding \emph{negative} formal completion of~$\frak m$,  defined analogously to $\widehat{\frak m}$ as above.


\section{The complete pro-unipotent group}\label{prounipgp}

The Lie algebra $\widehat{\mathfrak{m}}$ has an induced topology as a subset of 
$\mathfrak{X}:=\prod_{k\in\Z} \frak m_k$ with the product topology. 
The group $\Aut(\widehat{\frak m})\subseteq \End(\mathfrak{X})$ inherits a natural topology from $\widehat{\frak m}$ as a subspace of the space $\End(\mathfrak{X})$ of linear maps  with the pointwise topology.

For $y\in \widehat{\frak m}$, write the ${\frak m}_k$ component of $y$ as  $y_k $ and let $N$ be the smallest integer with $y_{N}\ne0$. Then  $y\in \widehat{\frak m}$ has the expression
$$y=\sum_{k=N}^\infty y_k.$$

For $k\in\mathbb{Z}$, set
\begin{align*}
\widehat{\mathfrak{m}}_{k}:= \begin{cases} 
      \widehat{\mathfrak{n}}_k&\text{if}\; k\ge 1 \\
     \bigoplus\limits_{{k\le j\le 0}} \mathfrak{m}_{j} \oplus \widehat{\frak n}^+ &\text{if}\; k\leq 0. \\
      \end{cases}
\end{align*}
where $\mathfrak{m}_{0}=\frak{gl}_2(-1)$. For $n\in\N$, we define
$$
  \widehat{{U}}_n := \left\{ \varphi\in\Aut(\widehat{\frak m}) \mid 
  \varphi(y)\in y + \widehat{\frak m}_{k+n}\text{ whenever }y\in{\frak m}_k,\;\text{for some}\; k\in\mathbb{Z}\right\}$$
and let 
$$\widehat{{U}}:=U_{-1}\widehat{{U}}_1,$$
where
$U_{-1}=\{ \exp(t\ad({e_{-1}})) \mid t\in \C^\times\}$ is the root group corresponding to the real root $\alpha_{-1}$.


Then $\widehat{U}$ is a closed pro-unipotent subgroup of $\Aut(\widehat{\frak m})$ 
with $\widehat{U}\ge\widehat{U}_1\ge\widehat{U}_2\ge\cdots$
and each $\widehat{U}/\widehat{U}_i$ is a finite dimensional unipotent algebraic group of automorphisms of 
$\widehat{\mathfrak{n}}/\widehat{\mathfrak{n}}_i$.

Every element  $ g\in\widehat{U}$ can be shown to have the form 
$$g=\prod_{i=1}^\infty \exp(\ad(x_i))$$
for some $x_i\in\widehat{\mathfrak{n}}_i$ (\cite{CJM}).
Any element $y\in \widehat{\frak m}$ is an infinite (formal) sum but, for $g\in \widehat{U}$, $g\cdot y$ is finite when $y$ is restricted to 
$\frak n^- \oplus \frak h \oplus (\widehat{\mathfrak n}^+/\widehat{\mathfrak{n}}_i)$
for all $i>0$. Thus every element $g\in \widehat{U}$  is pro-summable when expanded as an infinite sum of endomorphisms.

In particular,  $\widehat{U}$ is generated  (as a topological group) by the summable series
\begin{align*}
\exp(u\,\ad(e_{-1}))&=1 +u\,\ad(e_{-1})+\dfrac{u^2}{2!}\,\ad(e_{-1})^2+\cdots
\intertext{and the pro-summable  series}
\exp(u\,\ad(e_{jk})) &=1 +u\,\ad(e_{jk})+\dfrac{u^2}{2!}\,\ad(e_{jk})^2+\cdots,
\end{align*}
for $(j,k)\in I^{\im}$, all of which are power series in $u$ with constant term $1$.

The group $\widehat{U}$ is the analog of a completion of the unipotent subgroup of the adjoint form of 
a Chevalley or Kac--Moody group (as in \cite{CLM}). For Kac--Moody algebras, similar complete pro-unipotent groups have been constructed by Kumar (\cite{Ku}) and Rousseau (\cite{Rou}).


Our construction gives rise to the following analog of the group adjoint representation.
Let
 $$\Ad_{\widehat{U}/ \widehat{U}_i}:\widehat{U}/ \widehat{U}_i
\to\Aut(\widehat{\mathfrak{n}}/\widehat{\mathfrak{n}}_i)$$
be the adjoint representation of the finite dimensional linear algebraic group $\widehat{U}/ \widehat{U}_i$ over $\C$ and let
$$\Exp_i: \widehat{\mathfrak{n}}/\widehat{\mathfrak{n}}_i\to \widehat{U}/ \widehat{U}_i$$
be the exponential map. Taking inverse limits allows us to define $\Ad$ and $\Exp$ for $\widehat{U}$ and $\widehat{\frak{n}}^+$ (as in Section 4.4.25 of \cite{Ku}).

For $x\in \widehat{\mathfrak{n}}$ and $g\in \widehat{U}$ we have
$$\Exp(\Ad(g))(x)=g\Exp(x) g^{-1}$$
where $\Exp:\widehat{\frak{n}}^+\to \widehat{U}$ is the unique map making the following diagram commute

\begin{center}
\begin{tikzpicture}
  \node (A) {$\widehat{\frak{n}}^+$};
  \node (B) [below=of A] {$\widehat{U}$};
  \node (C) [right=of A] {$ \widehat{\mathfrak{n}}/\widehat{\mathfrak{n}}_i$};
  \node (D) [right=of B] {$\widehat{U}/ \widehat{U}_i$};
  \draw[-stealth] (A)-- node[left] {\small $\Exp$} (B);
  \draw[-stealth] (B)-- node [below] {\small $  $} (D);
  \draw[-stealth] (A)-- node [above] {\small $  $} (C);
  \draw[-stealth] (C)-- node [right] {\small $\Exp_i$.} (D);
\end{tikzpicture}
\end{center}


\section{Generators and relations}\label{gensrelnsgp}

In the finite dimensional case, 
Steinberg  (\cite{St}) gave a defining presentation for (adjoint and simply connected) Chevalley groups over commutative rings.
Tits (\cite{Ti}) gave generators and relations for Kac--Moody groups, generalizing the Steinberg presentation.  
A suitable generalization of these methods gives  us a group  $G(\frak m)$ for $\frak m$ in terms of generators and relations.

To motivate our method for constructing the group $G(\frak m)$, let $G$ be an adjoint Kac--Moody group corresponding to a symmetrizable Kac--Moody algebra $\mathfrak{g}$. Then  all  elements of $G$ can be constructed as automorphisms of  $\mathfrak{g}$. 

As we have seen,  it is not possible construct all the analogous group elements  as automorphisms of $\frak m$. 
For example, we may construct automorphisms of  $\widehat{\mathfrak{m}}$ corresponding to positive roots of $\frak m$ but not negative roots (as in the construction of $\widehat{U}$). We  may also consider groups $\GL_2(-1)$ and  $\GL_2({\ell,j,k})$  (defined below) which are automorphisms of $\frak{gl}_2$ subalgebras of $\mathfrak{m}$, corresponding to the real root $\alpha_{-1}$ and imaginary roots $\alpha_{\ell,jk}$ respectively. Using  pairwise intersections of  these subalgebras  allows us to construct $G(\mathfrak{m})$ as an amalgam of all these groups. This gives rise to the following presentation. Define the set of symbols
$$\mathcal{X} = \{H_1(s), H_2(s), X_{-1}(u),  Y_{-1}(u), X_{\ell,jk}(u), Y_{\ell,jk}(u) \mid s\in \C^\times, u\in \C, (\ell,j,k)\in E \}.$$
Define the constant
$$c_{\ell j} := (-1)^{\ell+1} \binom{j-1}{\ell} \, (\ell+1)(j-\ell)$$
and additional symbols:
\begin{align*}
 \widetilde{w}_{-1}(s) &:= X_{-1}(s) Y_{-1}(-s^{-1}) X_{-1}(s), \qquad \widetilde{w}_{-1} := \widetilde{w}_{-1}(1),\\
 \widetilde{w}_{\ell,jk}(s) &:= X_{\ell,jk}\left({s }\right) Y_{\ell,jk}\left(\frac{-s^{-1}}{{c_{\ell j}}}\right) 
  X_{\ell,jk}\left({s}\right),\qquad  \widetilde{w}_{\ell,jk} := \widetilde{w}_{\ell,jk}(1).
  \end{align*}
 Let $(g,h)=ghg^{-1}h^{-1}$ denote the usual group commutator.  We define the set of relations $\mathcal{R}$ as follows, for all $s,t\in \C^\times$, $u,v\in \C$, $(\ell,j,k),(m,p,q) \in E$.
\begin{align*}
\intertext{Relations in the $\GL_2({\C})$-subgroup associated with the real simple root:}
 X_{-1}(u)X_{-1}(v)&=X_{-1}(u+v),\\
Y_{-1}(u)Y_{-1}(v)&=Y_{-1}(u+v),\\
H_1(s)H_1(t)&=H_1(st),\\
H_2(s)H_2(t)&=H_2(st),\\
    H_1(s)H_2(t)&=H_2(t)H_1(s),\\
  \widetilde{w}_{-1}X_{-1}(u)\widetilde{w}_{-1}^{-1}&=Y_{-1}(-u),\\
  \widetilde{w}_{-1}Y_{-1}(u)\widetilde{w}_{-1}^{-1}&=X_{-1}(-u),\\
  Y_{-1}(-t)X_{-1}(s)Y_{-1}(t)&= X_{-1}(-t^{-1}) Y_{-1}(-t^{2}s)X_{-1}(t^{-1}),\\
  \widetilde{w}_{-1}(s)\widetilde{w}_{-1}&=H_1(-s)H_2(-s^{-1}), \\
  {\widetilde{w}_{-1}} H_1(s)\widetilde{w}_{-1}^{-1} &= H_2(s),\\
  {\widetilde{w}_{-1}} H_2(s)\widetilde{w}_{-1}^{-1} &= H_1(s),\\
  H_1(s)X_{-1}(u)H_1(s)^{-1}&=X_{-1}(su),\\
  H_2(s)X_{-1}(u)H_2(s)^{-1}&=X_{-1}(s^{-1}u),\\
  H_1(s)Y_{-1}(u)H_1(s)^{-1}&=Y_{-1}(s^{-1}u),\\
  H_2(s)Y_{-1}(u)H_2(s)^{-1}&=Y_{-1}(su),\\
\intertext{Relations between generators $X_{-1}(u),  Y_{-1}(u), X_{\ell,jk}(u), Y_{\ell,jk}(u)$:}
  (X_{\ell,jk}(u), Y_{m,pq}(v))&=1\qquad\text{for $j\neq p$, $k\neq q$, or $|\ell-m|>1$,}\\
  X_{\ell,jk}(u+v)&=X_{\ell,jk}(u)X_{\ell,jk}(v),\\
  Y_{\ell,jk}(u+v)&=Y_{\ell,jk}(u)Y_{\ell,jk}(v),\\
  (X_{-1}(s), X_{j-1,jk}(t))&=1,\\
  (Y_{-1}(s), X_{0,jk}(t))&=1,\\
  (X_{-1}(s), Y_{0,jk}(t))&=1,\\
  (Y_{-1}(s), Y_{j-1,jk}(t))&=1.\\
\intertext{Action of the element  $ \widetilde{w}_{-1}$ corresponding to the  real root $\a_{-1}$:}
    \widetilde{w}_{-1}X_{\ell,jk}(u)\widetilde{w}_{-1}^{-1}&=X_{j-1-\ell,jk}((-1)^{j-\ell-1} u), \\
  \widetilde{w}_{-1}Y_{\ell,jk}(u)\widetilde{w}_{-1}^{-1}&=Y_{j-1-\ell,jk}((-1)^{j-\ell-1} u), \\
\intertext{Relations in the $\GL_2(\C)$-subgroup associated with imaginary roots:}
  H_1(s) X_{\ell,jk}(u) H_1(s)^{-1}&=X_{\ell,jk}(s^{\ell+1}u),\\
  H_1(s) Y_{\ell,jk}(u) H_1(s)^{-1}&=Y_{\ell,jk}(s^{-(\ell+1)}u),\\
  H_2(s) X_{\ell,jk}(u) H_2(s)^{-1}&=X_{\ell,jk}(s^{j-\ell}u),\\
  H_2(s) Y_{\ell,jk}(u) H_2(s)^{-1}&=Y_{\ell,jk}(s^{-(j-\ell)}u),\\
  \widetilde{w}_{\ell,jk}(s)\widetilde{w}_{\ell,jk}&=H_1\left((-s)^{1/(\ell+1)}\right) H_2\left((-s)^{1/(j-\ell)}\right), \\
  Y_{\ell,jk}(-t)X_{\ell,jk}(s)Y_{\ell,jk}(t)&= 
  X_{\ell,jk}\left(\frac{-t^{-1}}{c_{\ell j}}\right)Y_{\ell,jk}\left(-c_{\ell j}t^2s\right)X_{\ell,jk}\left(\frac{t^{-1}}{c_{\ell j}}\right),\\
  \widetilde{w}_{\ell,jk} H_1(s) {\widetilde{w}_{\ell,jk}}^{-1}&=  H_2\left(s^{-(\ell+1)/(j-\ell)}\right),\\
  \widetilde{w}_{\ell,jk} H_2(s) {\widetilde{w}_{\ell,jk}}^{-1}&=  H_1\left(s^{-(j-\ell)/(\ell+1)}\right),\\
  \widetilde{w}_{\ell,jk} X_{\ell,jk}(u) {\widetilde{w}_{\ell,jk}}^{-1}&= Y_{\ell,jk}\left(\frac{-u}{c_{\ell j}}\right),\\
  \widetilde{w}_{\ell,jk} Y_{\ell,jk}(u) {\widetilde{w}_{\ell,jk}}^{-1}&= X_{\ell,jk}\left(-c_{\ell j}u\right),
\end{align*}
\begin{align*} 
Y_{\ell,jk}(s) =
X_{\ell,jk}\left(\frac{ s^{-1}}{c_{\ell j}}\right)
H_1\left(\left[-{{c_{\ell j}}}s\right]^{-1/(\ell+1)}\right)
H_2\left(\left[-{{c_{\ell j}}}s\right]^{-1/(j-\ell)}\right)\widetilde{w}_{\ell,jk}X_{\ell,jk}\left(\frac{ s^{-1}}{c_{\ell j}}\right).
\end{align*}

We now define $G(\mathfrak{m})$ as the group given by this presentation, that is,
$${G}(\mathfrak{m})=\langle \mathcal{X}\mid \mathcal{R}\rangle={F}(\mathcal{X})/N_\mathcal{R}$$
where ${F}(\mathcal{X})$ denotes the free group on $\mathcal X$ and $N_\mathcal{R}$ denotes the normal closure of the relations $\mathcal R$.

We define the following subgroups of $G(\mathfrak{m})$:
\begin{align*}
U^+(\C)&=\langle X_{-1}(u),\ X_{\ell,jk}(u) \mid  u\in\C, (\ell,j,k)\in E \rangle,\\
 U^+_{\im}(\C)&=\langle X_{\ell,jk}(u) \mid  u\in\C, (\ell,j,k)\in E \rangle,\\
U^+_{\im}(\Z)&=\langle X_{\ell,jk}(1)\mid (\ell,j,k)\in E\rangle,\\
H&=\langle H_1(s), H_2(s) \mid s\in \C^\times\rangle.
\end{align*}
We call $H$ the {\it toral} subgroup of  $G(\mathfrak{m})$. 


\begin{proposition}(\cite{CJM}) The subgroup $U^+_{\im}(\C)$  is a free product of additive abelian groups isomorphic to $\C$,  indexed over $E$. The subgroup $U^+_{\im}(\Z)$ is a countably-generated free group. 
For fixed $(\ell,j,k)$, the group \linebreak $\langle X_{\ell,jk}(u) \mid  a\in \Z\rangle\cong\Z$ is an infinite cyclic subgroup.
\end{proposition}

The following theorem gives the relationship between the groups $G(\frak m)$ and  $\widehat{U}$.
\begin{theorem}\label{main} The map 
$$ X_{\ell,jk}(u)\mapsto \exp(u\,\ad(e_{\ell,jk}))$$
embeds $U_{\im}^+(\C) $ as a dense subgroup
of
 $\widehat{U}$ . 
 \end{theorem}
Theorem~\ref{main}  shows that an element of the subgroup $U_{\im}^+(\C)$ of $G(\frak m)$ can be identified with an automorphism of $\widehat{\mathfrak{m}}$, and every
automorphism in $\widehat{U}$ can be approximated by elements of $U_{\im}^+(\C)$.

Our group $G(\mathfrak m)$ does not act as an  automorphism group of $\mathfrak{m}$.
However,  the action of the $\mathfrak{gl}_2(-1)$  subalgebra on $\mathfrak{m}$ is locally nilpotent, and 
thus elements of the 
$\GL_2(\C)$ subgroup
\begin{align*}\GL_2({-1})&=\langle X_{-1}(u), Y_{-1}(u), H_1(s), H_2(s)\mid u\in\C,\, s\in\C^\times\rangle
\end{align*}
can be identified with automorphisms of $\mathfrak{m}$. The groups $$\GL_2({\ell,j,k}) = \langle   X_{\ell,jk}(u),  Y_{\ell,jk}(u),H_1(s), H_2(s) \mid u\in\C,\, s\in\C^\times, (\ell,j,k)\in E \rangle$$ act as automorphisms on subalgebras $\frak{gl}_2(\ell,j,k)$ but not on all of $\frak m$.



\section{Analog of a simply connected Kac--Moody Chevalley group for $\frak m$}\label{simplyconn}

An analog of a simply connected  Chevalley group can be constructed for Kac--Moody algebras $\mathfrak g$ over commutative rings  (\cite{CG}, \cite{CLM}, \cite{Rou}, \cite{Ti}). Constructing these groups requires a significant amount of additional  data such as a $\mathbb{Z}$-form of the universal enveloping algebra of   $\mathfrak g$, as well as an integrable highest weight representation of $\mathfrak g$.

To construct an analog of a simply connected  Kac--Moody Chevalley group for $\frak m$  would require the use of an integrable highest weight representation of $\mathfrak m$. Recall that the simple imaginary root vectors which generate  $\frak m$ are not locally nilpotent on non-trivial highest weight $\frak m$-modules. 

In this section, we outline the construction of a group (\cite{CGJM}) associated to the parabolic subalgebra $\frak{gl}_2(-1)\oplus\frak u^+$ of $\frak m$ which acts locally nilpotently on a certain tensor algebra $T(\mathcal V)$, constructed in \cite{JLW}. The module $T(\mathcal V)$ is analogous to  standard irreducible modules  for semisimple and Kac--Moody algebras. 

Using the notation of Section~\ref{gensrelns}, for $j>0$, we define  
\begin{align*}
\mathcal{V}^-_j&=\coprod_{1\leq k\leq c(j)}{\mathcal{U}(\C f_{-1})\cdot f_{jk}}\\
\mathcal{V}^+_j&=\coprod_{1\leq k\leq c(j)}{\mathcal{U}(\C e_{-1})\cdot e_{jk}},
\end{align*}
where $\mathcal{U}$ denotes the universal enveloping algebra.  Set 
\begin{align*}
\mathcal{V}^-=\coprod_{j>0}{\mathcal{V}^-_j},\qquad \mathcal{V}^+=\coprod_{j>0}{\mathcal{V}^+_j}.
\end{align*}
Then we have $\mathfrak{u}^-=L(\mathcal{V}^-),$ $\mathfrak{u}^+=L(\mathcal{V}^+),$ where $L(V)$ denotes the free Lie algebra on a vector space $V.$
Let $\mathcal{V}$ denote $\mathcal{V}^-$.
 Let $T(\mathcal V)$ denote the tensor algebra of $\mathcal V$.

In \cite{JLW}, the authors showed that the Lie algebra $\mathfrak{m}$ can be realized as a Lie algebra of operators on the irreducible $\mathfrak{m}$-module $T(\mathcal{V}),$ identified  with a generalized Verma module induced from a one dimensional $\mathfrak{gl}_2$-representation $L(\lambda)$ for any $\lambda\in (\mathfrak{h}_A)^*$ satisfying 
	\begin{equation}\label{lambda-conditions}
	\lambda(h_{-1})=0,\quad \lambda(h_{jk})=a j +b
	\end{equation}
	for each $j>0,$ $1\leq k\leq c(j),$ for some fixed real numbers $a,b$ such that $a>0$ and $a+b>0.$


This representation has the property that the  parabolic  subalgebra $\frak{gl}_2(-1)\oplus\frak u^+$ of $\frak m$ is locally nilpotent on the $\frak m$-module $T(\mathcal V)$ (\cite{JLW}).


The module $T(\mathcal{V})$ has a filtration 
	\begin{align*}
	\C=T_0(\mathcal{V})\subseteq T_1(\mathcal{V})\subseteq T_2(\mathcal{V}) \subseteq \cdots \subseteq T_n(\mathcal{V})\subseteq \cdots \subseteq T(\mathcal{V}),
	\end{align*} 
	where $T_n(\mathcal{V})=T\left(\coprod_{j\leq n}\mathcal{V}_j\right)$, that is invariant under the action of $\mathfrak{gl}_2(-1)\oplus\frak u^+.$ In particular, every element of $T(\mathcal{V})$ lies in a proper $\mathfrak{gl}_2(-1)\oplus\frak u^+$-submodule that is a finitely generated tensor algebra.

Combining the degree grading on the tensor algebra $T(V)$ of any vector space $V$ given by 
\begin{align*}
T(V)&=\bigoplus_{s\geq 0}T^{(s)}(V) \text{ where}\\
T^{(s)}(V)&=\text{span}_\C\{v_1\otimes v_2\otimes \cdots \otimes v_s\mid v_1,\ldots,v_s\in V\},
\end{align*}
with the filtration on $T(\mathcal{V}),$ for any $n\geq 1$ we get the grading 
\begin{align*}
T_n(\mathcal{V})=\bigoplus_{s\geq 0}T^{(s)}_n(\mathcal{V}) \text{ where}\\
T^{(s)}_n(\mathcal{V})=T_n(\mathcal{V})\cap T^{(s)}(\mathcal{V}).
\end{align*}

The subspaces $T^{(s)}_n(\mathcal{V})$ are finite dimensional. The associated filtration on $T_n(\mathcal{V})$ is
\begin{align*}
&T^{(\leq 0)}_n(\mathcal{V})\subseteq T^{(\leq 1)}_n(\mathcal{V})\subseteq \cdots \subseteq T^{(\leq s)}_n(\mathcal{V})\subseteq \cdots\subseteq T_n(\mathcal{V}), \text{ where}\\
&T^{(\leq s)}_n(\mathcal{V}):=\coprod_{i=0}^{s}T^{(i)}_n(\mathcal{V}).
\end{align*}
We let $\rho_{n,s}=\rho|_{T^{(\leq s)}_n(\mathcal{V})}.$ For $n=0,$ define $T^{(0)}_0(\mathcal{V}):=\C$ and $T^{(i)}_0(\mathcal{V}):=0$ for $i>0.$ 

\begin{proposition}
	For any $n,s\geq 0$, the subspace $T^{(\leq s)}_n(\mathcal{V})$ is a $\mathfrak{gl}_2(-1)\oplus\frak u^+$-submodule, such that for $s\geq 1,$
	\begin{align*}
	\rho_{n,s}(\mathfrak{u}^+)T^{(s)}_n(\mathcal{V})&\subseteq T^{(s-1)}_n(\mathcal{V}),\\
	\rho_{n,s}(\mathfrak{gl}_2(-1))T^{(s)}_n(\mathcal{V})&\subseteq T^{(s)}_n(\mathcal{V}).
	\end{align*}
	In particular, $\rho_{n,s}(\mathfrak{u}^+)$ is  a finite dimensional nilpotent Lie algebra and $\rho_{n,s}(\mathfrak{gl}_2(-1))$ is finite dimensional.
\end{proposition}

Choosing  $\lambda$ subject to the conditions (\ref{lambda-conditions}) as above, for all $n,s\geq 0$, we may define a family of groups 
$G^\lambda_{n,s}(\mathfrak{gl}_2(-1)\oplus\frak u^+)$ associated to $\rho_{n,s}(\mathfrak{gl}_2(-1)\oplus\frak u^+)$ as follows (\cite{CGJM})
	\begin{align*}
		G^\lambda_{n,s}(\mathfrak{gl}_2(-1)\oplus\frak u^+)&=\langle \exp(c\rho_{n,s}(e_{jk})),\exp(c\rho_{n,s}(e_{-1})),\exp(c\rho_{n,s}(f_{-1})),\\
		&\exp(t\rho_{n,s}(h_{jk})) \exp(t\rho_{n,s}(h_1))\mid (j,k)\in I,\ c\in\C,\ t\in\C^\times \rangle.
	\end{align*}

	Since each $\rho_{n,s}(\mathfrak{u}^+)$ is finite dimensional, $G^\lambda_{n,s}(\mathfrak{gl}_2(-1)\oplus\frak u^+)$ is a well defined subgroup of $\text{Aut}(T_n^{\leq s}(\mathcal{V}))$. 
	This is an analog of  the parabolic subgroup of a simply connected  Chevalley group  associated to the subalgebra  $\frak{gl}_2(-1)\oplus\frak u^+$ of $\frak m$.

	We define subgroups of $G^\lambda_{n,s}(\mathfrak{gl}_2(-1)\oplus\frak u^+)$ as follows. Let $\mathfrak{n}^+=\mathfrak{u}^+\oplus \C e_{-1}$ as in Section~\ref{gensrelns} and define	\begin{align*}
U^\lambda_{n,s}(\mathfrak{n}^+)&=\langle \exp(c\rho_{n,s}(e_{jk})),\exp(c\rho_{n,s}(e_{-1})) \mid (j,k)\in I,\ c\in\C\rangle,\\
U^\lambda_{n,s}(\mathfrak{u}^+)&=\langle \exp(c\rho_{n,s}(e_{jk})) \mid (j,k)\in I,\ c\in\C\rangle,\\
H^\lambda_{n,s}&=\langle \exp(t\rho_{n,s}(h_1)), \exp(t\rho_{n,s}(h_2)) \mid t\in\C^\times\rangle,\\
L^\lambda_{n,s}&=\langle\exp(c\rho_{n,s}(e_{-1})),\exp(c\rho_{n,s}(f_{-1})), \exp(t\rho_{n,s}(h_1)), \exp(t\rho_{n,s}(h_2))\mid c\in\C,\ t\in\C^\times\rangle.
\end{align*}
	
\begin{theorem} (\cite{CGJM})  For $\lambda$ subject to conditions (\ref{lambda-conditions}) and for all $n,k\geq 0$ we have the following.
\begin{itemize}
\item The groups $U^\lambda_{n,s}(\mathfrak{n}^+)$ are pro-unipotent.
\item The groups $U^\lambda_{n,s}(\mathfrak{u}^+)$ are free products of additive groups isomorphic to $(\C,+)$ and are pro-unipotent.
\item The groups $H^\lambda_{n,s}$ are isomorphic to the toral subgroup of $\GL_2(\C)$.
\item The groups $L^\lambda_{n,s}$ are isomorphic to  $\GL_2(\C)$.
\item $G^\lambda_{n,s}(\mathfrak{gl}_2(-1)\oplus\frak u^+)=L^\lambda_{n,s} \ltimes U^\lambda_{n,s}(\mathfrak{n}^+)$.
	 \end{itemize}
	\end{theorem}


 \section{Conflicts of interest}
 On behalf of all authors, the corresponding author states that there is no conflict of interest. 
 \section{Data availability}
 This manuscript has no associated data.


\bigskip
\begin{thebibliography}{[CCCCC]}


\bibitem[ACMY]{ACMY} Ali, A., Carbone, L.,  Murray, S. H. and Yap, C. {\it Integral structure in the automorphism group of the Monster Lie algebra}, In preparation (2022) 

\bibitem[ACJM]{ACJM} Ali, A., Carbone, L., Jurisich, E. and Murray, S. H.,  {\it A complete group for the Borel subalgebra of a Borcherds algebras}, In preparation (2022)



\bibitem[B1]{B1} Borcherds, R. E., {\it Introduction to the Monster Lie algebra},  Groups, combinatorics and  geometry (Durham, 1990), 99--107, London Math. Soc. Lecture Note Ser., 165, Cambridge Univ. Press, Cambridge, (1992)

\bibitem[B2]{B2}  Borcherds, R. E., {\it Monstrous moonshine and monstrous Lie superalgebras},  Invent. Math. 109 (1992), no. 2, 405--444.


\bibitem[B3]{B3} Borcherds, R.E. {\it Generalized Kac-Moody algebras},  J. Algebra 115, 501--512 (1988)


\bibitem[Ca]{Ca} Carnahan, S., {\it Fricke Lie algebras and the genus zero property in Moonshine},  Journal of Physics A 50:40 (2017) (ArXiv).

\bibitem[CLM]{CLM} Carbone, L.,  Liu, D. and Murray, S. H., {\it Infinite dimensional Chevalley groups and Kac--Moody groups over $\Z$}, arXiv:1803.11204 [math.RT]

\bibitem[CCJMP]{CCJMP} Carbone, L., Chen, H., Jurisich, E., Murray, S. H. and Paquette, N.  {\it Imaginary reflections and automorphisms of the Monster Lie algebra}, In preparation (2022) 

\bibitem[CG]{CG} Carbone, L. and Garland, H., 
  {\it Existence of lattices in  Kac--Moody groups over finite fields},  Communications in Contemporary
Math, Vol 5, No.5,
 813--867,
(2003)

\bibitem[CGJM]{CGJM} Carbone, L., Ginory, A., Jurisich, E. and Murray, S. H., {\it Automorphism groups of standard induced modules for the Monster Lie Algebra}, In preparation  (2022)





\bibitem[CJM]{CJM} Carbone, L., Jurisich, E. and Murray, S. H.,  {\it A Lie group analog for the Monster Lie algebra}, In preparation (2022) 

 

\bibitem[CP]{CP}  Carbone, L. and Paquette, N. {\it Imaginary reflections and discrete symmetries in the heterotic Monster}, (2022) arXiv:2202.09720v3 [hep-th]


\bibitem[CN]{CN} J. H. Conway and S. P. Norton, {\it Monstrous Moonshine}, Bull. London Math. Soc. 11 (1979), 308--339.


\bibitem [Ga1]{Ga1} Garland, H.,  {\it{The arithmetic theory of loop groups}},
 Inst. Hautes tudes Sci. Publ. Math. No. 52 (1980), 5--136.
 
 \bibitem[GW]{GW}  Goodman, R. and Wallach, N. R. {\it Structure and unitary cocycle representations of loop groups and the group of diffeomorphisms of the circle},  J. Reine Angew. Math.  347  (1984), 69--133
 
\bibitem[HM]{HM} Harvey, J. and A. Moore, G., {\it On the algebras of BPS states},
Comm. Math. Phys. 197 (1998), no. 3, 489--519.

\bibitem[FLM1]{FLM1} Frenkel, Igor; Lepowsky, James; Meurman, Arne {\it Vertex operator algebras and the Monster}, Pure and Applied Mathematics, 134. Academic Press, Inc., Boston, MA, 1988


\bibitem[FLM2]{FLM2} Frenkel, I. B.; Lepowsky, J.; Meurman, A. A natural representation of the Fischer-Griess Monster with the modular function J as character. Proc. Nat. Acad. Sci. U.S.A. 81 (1984), no. 10, , Phys. Sci., 3256?3260



\bibitem[Jur1]{Jur1} Jurisich, E., {\it Generalized Kac-Moody Lie algebras, free Lie algebras and the structure of the Monster Lie algebra}, J. Pure Appl. Algebra 126 (1998), no. 1-3, 233--266

\bibitem[Jur2]{Jur2} Jurisich, E.,  {\it An exposition of generalized Kac-Moody algebras. Lie algebras and their representations}, (Seoul, 1995), 121--159, Contemp. Math., 194, Amer. Math. Soc., Providence, RI, (1996).

 \bibitem[JLW]{JLW} Jurisich, E., Lepowsky, J. and Wilson, R. L., {\it Realizations of the Monster Lie algebra},  Selecta Math. (N.S.) 1 (1995), no. 1, 129--161
 
   \bibitem[KP]{KP}   Kac, V. and Peterson, D., {\it  Defining relations of certain infinite dimensional groups}, Ast\'erisque Num\'ero Hors S\'erie (1985), 165--208. 
   
\bibitem[Ku]{Ku} Kumar, S., {\it Kac--Moody groups, their flag varieties and representation theory}, Vol. 204, Progress in Mathematics, Springer,
(2012).


\bibitem[LL]{LL} Lepowsky, J.and  Li, H. {\it Introduction to vertex operator algebras and their representations}, Progress in Mathematics, 227. Birkha\"user Boston, Inc., Boston, MA, 2004. xiv+318 pp

\bibitem[Ma]{Ma}  Marcuson, R. {\it Tits' systems in generalized nonadjoint Chevalley groups}  J. Algebra  34  (1975), 84--96



 \bibitem[Mar]{Mar} Marquis, T.  {\it An introduction to Kac-Moody groups over fields},  EMS Textbooks in Mathematics. European Mathematical Society (EMS), Z\"urich, 2018. xi+331 pp. 

\bibitem[Ma84]{Ma84} Mathieu, O. {\it Sur la construction de groupes associ\'es aux alg\'ebre
s de Kac--Moody}. (French) [On the construction of groups associated with Kac--Moody algebras]  C. R. Acad. Sci. Paris S\'er. I Math.  299  (1984),  no. 6, 161--164

\bibitem[MT]{MT}  Moody, R. V. and Teo, K. L. {\it Tits' systems with crystallographic Weyl groups}, J. Algebra 21 (1972), 178--190

\bibitem[PPV]{PPV} Paquette, N. M., Persson, D. and  Volpato, R., {\it  BPS algebras, genus zero and the heterotic Monster}, J. Phys. A 50 (2017), no. 41, 414001, 17 pp

\bibitem[Rou]{Rou} Rousseau, G., {\it Groupes de Kac-Moody d\'eploy\'s sur un corps local II. Masures ordonn\'es. (French) [Split Kac-Moody groups over a local field II. Ordered hovels]}, Bull. Soc. Math. France 144 (2016), no. 4, 613--692.




 \bibitem[Sl]{Sl} Slodowy, P. {\it An adjoint quotient for certain groups attached to Kac--Moody algebras},  infinite dimensional groups with applications (Berkeley, Calif., 1984),  307--333, Math. Sci. Res. Inst. Publ., 4, Springer, New York, (1985)
 
 
 

\bibitem[St]{St} Steinberg, R., {\it G\'en\'erateurs, relations et rev\^etements de groupes alg\'ebriques}, (French) 1962 Colloq. Th\'eorie des Groupes Alg\'ebriques (Bruxelles, 1962) pp. 113--127 Librairie Universitaire, Louvain; Gauthier-Villars, Paris
 
  \bibitem[Ti]{Ti}   Tits, J., {\it Uniqueness and presentation of Kac--Moody groups
over fields}, Journal of Algebra 105, 542--573 (1987)

\end{thebibliography}
\end{document}